\documentclass[12pt]{article}

\usepackage{amsmath,amsthm,amsfonts,latexsym,amsopn,verbatim,amscd,amssymb}
\newenvironment{psmallmatrix}
  {\left(\begin{smallmatrix}}
  {\end{smallmatrix}\right)}

\theoremstyle{plain}
\newtheorem{Thm}{Theorem}[section]
\newtheorem{Lem}[Thm]{Lemma}

\theoremstyle{definition}

\numberwithin{equation}{section}

\DeclareMathOperator{\idem}{idem}
\DeclareMathOperator{\nil}{nil}
\DeclareMathOperator{\U}{U}

\DeclareMathOperator{\Nin}{Nin}

\newcommand{\bnum}{\begin{enumerate}}
\newcommand{\enum}{\end{enumerate}}
\begin{document}
\begin{center}
\textbf{Formal triangular matrix ring with nil clean index $4$ }\\
\end{center}
\begin{center}
Dhiren Kumar Basnet\\
\small{\it Department of Mathematical Sciences, Tezpur University,
 \\ Napaam, Tezpur-784028, Assam, India.\\
Email: dbasnet@tezu.ernet.in}
\end{center}
\begin{center}
Jayanta Bhattacharyya \\
\small{\it Department of Mathematical Sciences, Tezpur University,
 \\ Napaam, Tezpur-784028, Assam, India.\\
Email: jbhatta@tezu.ernet.in}
\end{center}
\noindent \textit{\small{\textbf{Abstract:}  }} For an element $a \in R$, let $\eta(a)=\{e\in R\mid e^2=e\mbox{ and }a-e\in \mbox{nil}(R)\}$. The nil clean index of $R$, denoted by $\Nin(R)$, is defined as $\Nin(R)=\sup \{\mid \eta(a)\mid: a\in R\}$. In this article we have characterized formal triangular matrix ring $\begin{psmallmatrix}A & M\\0 & B\end{psmallmatrix}$ with nil clean index $4$.
\bigskip

\smallskip

\noindent \small{\textbf{\textit{$2010$ Mathematics Subject Classification:}} 16U99} \\
\smallskip

\bigskip
\section{Introduction}
Throughout this article $R$ denotes a associative ring with unity. The set of nilpotents and set of idempotents are denoted by $\nil(R)$ and $\idem(R)$ respectively. The cyclic group of order $n$ is denoted by $C_n$ and $|S|$ denotes the cardinality of the set $S$. For an element $a\in R$, if $a - e\in \nil(R)$ for some $e\in \idem(R)$, then $a = e +( a - e)$ is said to be a nil clean expression of $a$ in $R$ and $a$ is called a nil clean element\cite{nc,ajd}. The ring R is called nil clean if each of its elements is nil clean.

\par For an element $a \in R$, let $\eta(a)=\{e\in R\mid e^2=e\mbox{ and }a-e\in \mbox{nil}(R)\}$. The nil clean index of $R$, denoted by $\Nin(R)$, is defined as $\Nin(R)=\sup \{\mid \eta(a)\mid: a\in R\}$\cite{nci}. Characterization of arbitrary ring with nil clean indices $1$, $2$ and few sufficient condition for a ring to be of nil clean index $3$ is given in \cite{nci}. In this article we have characterized formal triangular matrix ring $\begin{psmallmatrix}A & M\\0 & B\end{psmallmatrix}$ with nil clean index $4$, where $A$ and $B$ are rings and $M $ is a $A-B-$bimodule. Following results about nil clean index will be used in this article.

\begin{Lem}\label{la11}
([$1$, Lemma $2.5$]) Let $R=\left(\begin{array}{cc}  A & M \\  0 & B \\\end{array}\right)$, where $A$ and $B$ are rings, $_AM_B$ is a bimodule. Let $\mbox{Nin}(A)=n$ and $\mbox{Nin}(B)=m$. Then
\begin{enumerate}
  \item $\mbox{Nin}(R)\geq |M|.$
\item If $(M,+)\cong C_{p^k}$, where $p$ is a prime and $k\geq1$, then
$\mbox{Nin}(R)\geq n+[\frac{n}{2})(|M|-1)$, where $[\frac{n}{2})$ denotes the least integer greater than or equal to $\frac{n}{2}$.
\item Either $\mbox{Nin}(R)\geq nm+|M|-1$ or $\mbox{Nin}(R)\geq 2nm$.
\end{enumerate}
\end{Lem}

\begin{Lem}\label{la12}
([$1$, Lemma $2.6$]) Let $R=\left(\begin{array}{cc}  A & M \\  0 & B \\\end{array}\right)$, where $A$ and $B$ are rings, $_AM_B$ is a bimodule with $(M,+)\cong C_{2^r}$. Then $\mbox{Nin}(R)=2^r\mbox{Nin}(A)\mbox{Nin}(B)$.
\end{Lem}

\begin{Thm}\label{4}
([$1$, Theorem $4.1$]) $\mbox{Nin}(R)=2$ if and only if $R=\left(
                                    \begin{array}{cc}
                                      A & M \\
                                      0 & B \\
                                    \end{array}
                                  \right),$ where $\mbox{Nin} (A)=\mbox{Nin} (B)=1$ and $_AM_B$ is a bimodule with $|M|=2.$
\end{Thm}
\begin{Thm} \label{th1} 
([$1$, Proposition $4.2$]) If $R=\left(\begin{array}{cc}A & M \\0 & B \\\end{array} \right),$ where $\mbox{Nin}(A)=\mbox{Nin}(B)=1$ and $_AM_B$ is a bimodule with $|M|=3$ then $\mbox{Nin}(R)=3.$
\end{Thm}

\section{Main result}
\begin{Thm} Let $R=\left(\begin{array}{cc} A & M \\ 0 & B \\\end{array}\right)$, where $A$ and $B$ are rings, $_AM_B$ is a non trivial bimodule. Then $\Nin(R)=4$ if and only if one of the following holds:
\begin{enumerate}
  \item[$(1)$] $(M,+)\cong C_2$ and $\Nin(A)\Nin(B)=2$.
  \item[$(2)$] $(M,+)\cong C_4$ and $\Nin(A)=\Nin(B)=1$.
  \item[$(3)$] $(M,+)\cong C_2\oplus C_2$ plus one of the following
  \begin{enumerate}
  \item[$(a)$] $\Nin(A)=\Nin(B)=1$.

  \item[$(b)$] $\Nin(A)=1,~~
  B=\left(\begin{array}{cc} S & W \\ 0 & T \\\end{array}\right)$, where $\Nin(S)=\Nin(T)=1$ and $|W|=2$, and $eM(1_B-f)+ (1_A-e)Mf \neq 0$ for all $e^2=e\in A$ and $f\in\eta(b)$, where $b\in B$ with $|\eta(b)|=2$.

  \item[$(c)$] $\Nin(B)=1,~~
  A=\left(\begin{array}{cc} S & W \\ 0 & T \\\end{array}\right)$, where $\Nin(S)=\Nin(T)=1$ and $|W|=2$, and $eM(1_B-f)+ (1_A-e)Mf \neq 0$ for all $e^2=e\in B$ and $f\in\eta(a)$, where $a\in A$ with $|\eta(a)|=2$.
  \end{enumerate}
\end{enumerate}
\end{Thm}
\noindent$Proof:~~(\Leftarrow)$ If $(1)$ holds then by \textbf{Lemma \ref{la12}}, we get $\Nin(R)=4$.\\
If $(2)$ holds then $\Nin(R)\geq |M|=4$. Now, for any $\alpha=\left(\begin{array}{cc} a & x\\ 0 & b \\\end{array}\right)\in R$,

$$\eta(\alpha)=\left\{\left(\begin{array}{cc} e & w \\ 0 & f \\\end{array}\right)\in R:~ e\in\eta(a),~f\in\eta(b),~w=ew+fw\right\}.$$

\noindent Because $|M|=4,~|\eta(a)|\leq 1$ and $|\eta(b)|\leq 1$, it follows that
$|\eta(\alpha)|\leq 4$. Hence $\Nin(R)=4$.\\

\noindent Let $(3)~(a)$ holds, then $\Nin(R)\geq |M|=4$. Now, for any $\alpha=\left(\begin{array}{cc} a & x\\ 0 & b \\\end{array}\right)\in R$,
$$\eta(\alpha)=\left\{\left(\begin{array}{cc} e & w \\ 0 & f \\\end{array}\right)\in R:~ e\in\eta(a),~f\in\eta(b),~w=ew+fw\right\}.$$Because $|M|=4,~|\eta(a)|\leq 1$ and $|\eta(b)|\leq 1$, it follows that $|\eta(\alpha)|\leq 4$. Hence $\Nin(R)=4$.

\noindent Suppose $(3)(c)$ holds, then clearly $\Nin(R)\geq|M|=4.$ Let
$\alpha=\left(\begin{array}{cc} a & w \\ 0 & b \\\end{array}\right)\in R$. We show that $|\eta(\alpha)|\leq 4$ and hence $\Nin(R)=4$ holds. Since $\Nin(B)=1$, we can assume that $\eta(b)=\{f_0\}$. Then as above we have
$$\eta(\alpha)=\left\{\left(\begin{array}{cc} e & z \\ 0 & f_0 \\\end{array}\right)\in R: e\in\eta(a), z=ez+zf_0 \right\}$$

\noindent If $|\eta(a)|\leq1$, then $|\eta(\alpha)|\leq|\eta(a)|.|M|\leq4$. So we can assume that $|\eta(a)|=2$. Write $\eta(a)=\{e_1,~e_2\}$. Thus $\eta(\alpha)=T_1\bigcup T_2$, where
$$T_i=\left\{\left(\begin{array}{cc} e_i & z \\ 0 & f_0 \\\end{array}\right)\in R: (1_A-e_i)z=zf_0 \right\}\hspace{1.cm}(i=1,~2).$$
\noindent Since $\eta(1_A-a)=\{1_A-e_1,~1_A-e_2\}$, the assumption $(3)(c)$ shows that
$\{z\in M:(1_A-e_i)z=zf_0\}$ is a proper subgroup of $(M,+)$; so $|T_i|\leq2$ for $i=1,~2$. Hence $|\eta(\alpha)|\leq|T_1|+|T_2|\leq 4$

\noindent $(\Rightarrow)$ Suppose $\Nin(R)=4$. Then $2\leq|M|\leq \Nin(R)=4$. If $|M|=2$ then $\Nin(A)\Nin(B)=2$ by \textbf{Lemma \ref{la12}}, so $(1)$ holds.

\noindent Suppose $|M|=3$, then we have by \textbf{Lemma \ref{la11}}, $\Nin(A)+ |M|\leq \Nin(R)$, showing $\Nin(A)\leq 2$, similarly $\Nin(B)\leq 2$. But $\Nin(A)=2=\Nin(B)$ will give $\Nin(R)\geq6$ by \textbf{Lemma \ref{la11}} and $\Nin(A)=\Nin(B)=1$ will give $\Nin(R)=3$ by \textbf{Theorem \ref{th1}}. Hence the only possibility is $\Nin(A)\Nin(B)=2$, so without loss of generality we assume that $\Nin(A)=2$ and $\Nin(B)=1$. Write $M =\{ 0, x, 2x\}$. Now by \textbf{Theorem \ref{4}}, we have
$A=\left(\begin{array}{cc}T & N \\0 & S \\\end{array}\right)$, where $T$ \& $S$ are rings, $_TN_S$ is bimodule with $\Nin(T)=\Nin(S)=1$ and $|N|=2$. Note that for $e\in \idem(A),~ ex\in\{0,~x\}$, for if $ex=2x$, we have $2x=ex=e(ex)=e(2x)=e(x+x)=ex+ex=2x+2x=4x=x$ which is not true. 

\noindent Now Let $a=\left(\begin{array}{cc}1_T & 0 \\0 & 0 \\\end{array}\right)\in A$ such that\\
$a=\left(\begin{array}{cc}1_T & 0 \\0 & 0 \\\end{array}\right)+
\left(\begin{array}{cc}0 & 0 \\0 & 0 \\\end{array}\right)=
\left(\begin{array}{cc}1_T & y \\0 & 0 \\\end{array}\right)+
\left(\begin{array}{cc}0 & -y \\0 & 0 \\\end{array}\right)$.\\
\noindent Let us denote, $e_1=\left(\begin{array}{cc}1_T & 0 \\0 & 0 \\\end{array}\right),~
e_2=\left(\begin{array}{cc}1_T & y \\0 & 0 \\\end{array}\right),~
n_1=\left(\begin{array}{cc}0 & 0 \\0 & 0 \\\end{array}\right)$ and\\ $n_2=\left(\begin{array}{cc}0 & -y \\0 & 0 \\\end{array}\right)$, where $e_1,~e_2\in\idem(A)~\&~n_1,~n_2\in\nil(R)$. Now we have following cases:\\
\noindent \textbf{Case I:} Let $e_1x=e_2x=0$, then we have an element \\
$\beta= \left(\begin{array}{cc}(1_A-a) & 0 \\0 & 0 \\\end{array}\right)\in R$ such that
\begin{align}
\beta=&\left(\begin{array}{cc}(1_A-e_1) & z \\0 & 0 \\\end{array}\right)+
\left(\begin{array}{cc}-n_1 & -x \\0 & 0 \\\end{array}\right)\hspace{2cm} \forall z\in M \notag\\
=&\left(\begin{array}{cc}(1_A-e_2) & z \\0 & 0 \\\end{array}\right)+
\left(\begin{array}{cc}-n_2 & -x \\0 & 0 \\\end{array}\right)\hspace{2cm} \forall z\in M\notag
\end{align}
are six nil clean expressions for $\beta,$ which implies $|\eta(\beta)|\geq 6$, that is $\Nin(R)\geq6$, which is not possible.

\noindent \textbf{Case II:} Let $e_1x=e_2x=x$, then we have an element \\
$\alpha= \left(\begin{array}{cc}a & 0 \\0 & 0 \\\end{array}\right)\in R$ such that
\begin{align}
\alpha=&\left(\begin{array}{cc}e_1 & z \\0 & 0 \\\end{array}\right)+
\left(\begin{array}{cc}n_1 & -x \\0 & 0 \\\end{array}\right)\hspace{2cm} \forall z\in M \notag\\
=&\left(\begin{array}{cc}e_2 & z \\0 & 0 \\\end{array}\right)+
\left(\begin{array}{cc}n_2 & -x \\0 & 0 \\\end{array}\right)\hspace{2cm} \forall z\in M\notag
\end{align}
are six nil clean expressions for $\alpha,$ which implies $|\eta(\alpha)|\geq 6$, that is $\Nin(R)\geq6$, which is also not possible.\\

\noindent \textbf{Case III:} Let $e_1x=x$ and $e_2x=0$, then we have $(e_1-e_2)x=x$.\\
Let $j=e_1-e_2$, then clearly $j\in\nil(A)$ and we have \\
$jx=x~~\Rightarrow~~(1_A-j)x=0~~\Rightarrow~~x=0,~~($ as $(1-j)\in\U(A))$.\\
Which is not possible.\\

\noindent \textbf{Case IV:} Let $e_1x=0$ and $e_2x=x$, as in case III, we get a contradiction.\\
Hence if $M\cong C_3,~~\Nin(R)$ is never 3.

\noindent Suppose $|M|=4$. If $(M,+)\cong C_4$, then $\Nin(A)\Nin(B)=1$ by \\Lemma \ref{la12}, So $(2)$ holds. Let $(M,+)\cong C_2\oplus C_2$. Since $\Nin(R)=4$, by \textbf{Lemma \ref{la11}}, we have $\Nin(A)\Nin(B)\leq2$. If $\Nin(A)\Nin(B)=1$ then $(3)(a)$ holds. If $\Nin(A)\Nin(B)=2$, without loss of generality we can assume $\Nin(A)=2$ and $\Nin(B)=1$. So by \textbf{Theorem \ref{4}}, we have $A=\left(\begin{array}{cc}S & W \\0 & T \\\end{array}\right)$ where $\Nin(S)=\Nin(T)=1$, and $|W|=2$. To complete the proof suppose in contrary that $eM(1_B-f)+(1_A-e)Mf=0$ for some $f^2=f\in B$ and $e\in \eta(a)$, where $a\in A$ with $|\eta(a)|=2.$ Then $ew=wf$ for all $w\in M$. It is easy to check that $\eta(a)=\{e,~e+j\}$ where $j=\left(\begin{array}{cc}0 & w_0 \\0 & 0 \\\end{array}\right)\in A$ with $0\neq w_0\in W$. Thus, for $\gamma:=\left(\begin{array}{cc}1_A-e & 0 \\0 & f \\\end{array}\right)$,
$$\eta(\gamma)\supseteq \left\{\left(\begin{array}{cc}1_A-e & 0 \\0 & f \\\end{array}\right),\left(\begin{array}{cc}1_A-(e+j) & 0 \\0 & f \\\end{array}\right),\left(\begin{array}{cc}1_A-e & w \\0 & f \\\end{array}\right): w\in M\right\}$$ So $|\eta(\gamma)|\geq 5$, a contradiction. Hence $(3)(c)$ holds, similarly $(3)(b)$ can be proved$\hfill \Box$


\begin{thebibliography}{99}
\bibitem{nci}Basnet, D. Kr. and Bhattacharyya, J. Nil clean index of rings Internationa, \textit{I. Electronic Journal of Algebra}, 15, 145-156, 2014.
    
\bibitem{ajd}Diesl, A. J., \textit{Classes of Strongly Clean Rings}. Ph D. Thesis, University of California, Berkeley, 2006.
    
\bibitem{nc}Diesl, A. J., Nil clean rings .∗\textit{Journal of Algebra}, 383 : 197 - 211, 2013.

\bibitem{wkn} Nicholson, W.K., Lifting idempotents and exchange rings. \textit{Transactions of the
american mathematical society}, 229 : 269 - 278, 1977.

\end{thebibliography}
\end{document}